\documentstyle[11pt,amssymb,twoside]{article}
\font\elevengtc=eufm10 scaled \magstephalf
\font\ninegtc=eufm9

\font\sevengtc=eufm7

\newfam\gtcfam
\def\gtc{\fam\gtcfam\elevengtc}
\textfont\gtcfam=\elevengtc
\scriptfont\gtcfam=\ninegtc
\scriptscriptfont\gtcfam=\sevengtc
   \def\ZZ{{\Bbb Z}}
   \def\CC{{\Bbb C}}    \def\PP{{\Bbb P}}
\def\hom{\mathop{\rm Hom}\nolimits}
\def\End{\mathop{\rm End}\nolimits}
\def\ext{\mathop{\rm Ext}\nolimits}
\def\Pic{\mathop{\rm Pic}\nolimits}

\def\le{\leqslant}
\def\ge{\geqslant}
\title{Testing $S$-duality conjecture and exceptional bundles.}
\author{B.V. Karpov}
\date{\today}
\begin{document}
\thispagestyle{empty}
\maketitle

\begin{abstract}
In this paper we deal with algebro-geometrical problems connected with 
testing $S$-duality conjecture for super-symmetric Yang-Mills quantum field 
theories in four dimensions. We describe all field configurations such
that beta function coefficient for gauge groups of rank 1 and 2 is zero or
negative. 
We also give special series of such field configurations for gauge groups of an arbitrary 
rank. Realization of one of the series discovers a 
connection 
%
%
with exceptional bundles. That points again at relation between $S$-duality
and string duality.
\end{abstract}

{\large\bf Introduction.}

\vspace{2ex}

Testing $S$-duality conjecture for super-symmetric Yang-Mills quantum field 
theories in four dimensions is a geometric problem in a nature
(one works as a rule with the case $SYM$,$N=4$ with or without 
twisting, or $\ N=2\ +$ supermultiplet). 
Recall that for  possibility of calculation of quantum 
correlation functions in semi-classical limit
it is necessary for the value of the coefficient of 
Gell-Mann --- Low $\beta$-function 
in 
one-loop approximation
on
field configuration of the theory
to 
be 
zero (or at least negative). In this case the gauge coupling constant 
(electro-magnetic $e$ or strong $g$) is well-defined. Besides, if the
beta coefficient is zero, then there is no anomalous $U(1)$ global symmetry,
so the physics depends on a vacuum angle $\theta$ and our theory is
scaling-invariant. Therefore, there appears a constant as a 
parameter of our theory
$$ \tau=\frac{\theta}{2\pi}+\frac{4\pi i}{g^2} $$    
(belonging to the higher hyperplane) and the 
$S$-duality conjecture predicates that the 
partition function
$Z_\tau$ of our theory 
admits some modular properties under integral linear fractional 
transformations of $\tau$. The function $Z_\tau$ is a power series of
$q=e^{2\pi i\tau}$, so, really the question is about behavior of the
partition function under $\tau\mapsto -\frac 1\tau$.
Assumption about 
this symmetry (for $\theta=0$) was made by Montonen and 
Olive [MO] more than twenty years old. To check this pure mathematical 
fact we mimic Vafa and Witten's method [VW], which is applied 
unfortunately only to partial cases of 4-varieties (algebraic surfaces del 
Pezzo and K3) for which theorems of vanishing harmonic spinors in 
canonical and adjoint representations hold. We wouldn't use the vanishing
arguments proposed in [VW] for the general case because of the problems 
with compactification. We hope to consider it in the subsequent
articles.

The aim of this paper is to put mathematical problems that are related to
testing $S$-duality conjecture and to show that exceptional bundles 
naturally appear in this subject. 

In the first section we describe all configurations of fields
such that the coefficient of $\beta$-function is zero or negative
for small (1 and 2) ranks of the fiberwise
gauge group and point out some series of such field configurations
for an arbitrary 
rank of the fiberwise gauge group (the series a), b) and c) in item 3 of
{\bf 1.4.} 
In particular, we show that the complete program 
of testing interesting for physicists $S$-duality conjecture in "high
energies" for  the gauge groups of rank 1 and 2 comes to five cases
(in which mathematical research are being run now, see, for example,
[GL], [Ku1], [Ku2]). For an arbitrary rank 
of the gauge group three series are given with two of which physicists are
working (for the series a) see~[VW] and for b) see~[APS]).

In the second section we discuss SYM QFT over del Pezzo surfaces
with the field configurations a) and b) given in
Sec.~1. We use the interpretation of the spaces of solutions up to gauge 
of the 
field equations as a moduli spaces of stable bundles. Under this, the 
coefficients of the partition functions are interpreted as the degree
of the top Chern class of some standard bundles over the moduli spaces
(just like the top Chern class of the cotangent bundle sends us to the
topological Euler class in [VW]).

In the third section we consider the arithmetical conditions 
on the topological invariants of the stable bundles that appear in Sec.~2.
We show that the exceptional bundles with zero cohomology groups 
under some inequalities on the topological invariants are initial points
of sequences of sheaf moduli spaces that give partition functions
suitable for testing the $S$-duality conjecture.

In the last section we put an algebro-geometrical problem for c) series,
but without physical motivation.

It is necessary to observe that there is a parallel program of testing 
$S$-duality, in "low energies", connected with pure geometric 
description for dynamics of Coulomb branch of vacua ([DW], [APS]).
We don't discuss this approach here. In our cases non-abelian gauge symmetry
is unbroken. 

The author is grateful to Max-Planck-Institut 
f\"ur
Mathematik in Bonn for hospitality and 
to A.~N.~Tyurin for numerous useful discussions
without which the work on this paper was impossible.

The condition on the beta function coefficient gives the
main selection for field configurations in QFT that are suitable for testing
$S$-duality in high energies. We begin just with this.

\vspace{2ex}

{\large\bf 1. Beta function coef\/f\/icient.}

\nopagebreak[3]

\vspace{2ex}

Mathematical subject-matter of this section is certainly
well known for specialists in representation theory. We give it here only
for convenient transition to geometry in the sequel.

{\bf 1.1.} Let $\gamma$ be an irreducible finite dimensional complex 
representation of $su(r)$ (the Lie algebra for the fiberwise gauge group of our 
theory) and let $b_{\gamma}$ be the Casimir element of this 
representation. It is a central element of the universal enveloping algebra,
so, all $b_{\gamma}$ 
are proportional to each other. 

Let $can$ be the canonical representation of $su(r)$ in $\CC^{\,r}$.
Denote by $c_2(\gamma)$ the proportionality coefficient
between $b_{can}$ and $b_{\gamma}$:  
$$b_{can}= c_2(\gamma)b_{\gamma}\,. $$
Any complex representation of $su(r)$ can be extended to its 
complexification $sl(r,\CC)$ without changing $c_2$, so we consider 
$\gamma$ as an
irreducible finite dimensional representation of $sl(r,\CC)$. Let 
$\{e_i\}$ be a basis in $sl(r,\CC)$.
Consider the dual basis $\{f_i\}$ with respect to the bilinear form
\begin{equation}\label{bxy}
B_{\gamma} (X;Y)=Tr(\gamma (X)\,\gamma(Y))\,,
\end{equation}
i.e. such that
$B_{\gamma} (e_i;f_j)=\delta_{ij}$.
Recall that by definition, $b_{\gamma}=\sum\limits_i e_if_i$. 
All bilinear forms (\ref{bxy}) are also proportional to each other as 
invariant forms on simple Lie algebra. Therefore, we have:
\begin{equation}\label{c2viabxy}
B_{\gamma}=c_2(\gamma)B_{can}\,.
\end{equation}
For the dual representation  
$B_{\gamma^*} (X;Y)=B_{\gamma} (X;Y)$, so the following equality holds:
\begin{equation}\label{dual}
c_2(\gamma^*)=c_2(\gamma)\,.
\end{equation}

\vspace{1ex}

{\bf Remark.} One can define $c_2$ for an irreducible $u(r)$-representation
$\widetilde\gamma$ as $c_2(\gamma)$, where $\gamma$ is the restriction of
$\widetilde\gamma$ on the subalgebra $su(r)\subset u(r)$. Clearly, all the results
of this section hold also for $u(r)$-representations.

\vspace{1ex}

{\bf Remark.} Though notation $c_2(\gamma)$ is not quite conventional in 
the theory of Lie algebras, it justifies itself in subsequent transition 
to geometry. Namely, if the representation $\gamma$  is considered as a 
tensor operation on vector bundles of rank $r$, then for such a bundle $E$
$c_2(\gamma)$ is the coefficient 
at $c_2(E)$ in the expression of $c_2(\gamma(E))$ via the Chern classes of
$E$.

\vspace{1ex}

Consider a finite collection of representations
$\{\gamma_1,\dots,\gamma_{N_a}\}$,
among which there are some isomorphic to each other.
From physical viewpoint this collection describes configuration of fields
for QFT under consideration (see Sec.~2). 
The coefficient of beta function takes the collection 
$\{\gamma_1,\dots,\gamma_{N_a}\}$ to the value
$$
\beta(\{\gamma_i\})=-c_2(ad) + \sum_i c_2(\gamma_i)\,.
$$  

The aim of this section is to describe all cases when 
$\beta(\{\gamma_i\})\le 0$ for algebras $su(2)$ and $su(3)$ and to point 
out some series of such cases for $su(r)$.

\vspace{1ex}

{\bf Remark.} Since we deal with Lie algebras, freedom of choice of the gauge 
group is still remained for us. For instance, $su(2)\cong so(3)$,  that 
corresponds to choosing the determinant of the bundle modulo its rank.

\vspace{1ex}

There is  a general formula which expresses 
$c_2(\gamma)$ in terms of  the higher weight of the representation $\gamma$. 

\vskip 10pt

{\bf 1.2. Proposition.} {\it Let $R$ be the root system for $sl(r,\CC)$, 
$\lambda$ be the highest weight of $\gamma$ with respect to a base 
of $R$, $\rho$ be the half-sum of the positive roots, and
$(\ ;\ )$ be the Killing form. Then}
\begin{equation}\label{kasimir}
c_2(\gamma)=\frac{(\lambda;\lambda+2\rho)}{r^2-1}\dim\Gamma=
\frac{(\lambda;\lambda+2\rho)}{r^2-1}  
\prod_{\alpha\in R_+}\frac{(\lambda+\rho;\alpha)}{(\rho;\alpha)}\;.
\end{equation}
{\it Here $\Gamma$ is the representation space for $\gamma$.}

\vskip 10pt

{\sc The proof} is based on two formulas of representation theory of
semisimple Lie algebras. The former is the Freudenthal formula 
applied to the highest weight (see, for example, [FH, (25.14)]) which 
states that the Casimir element $b_{can}$ acts on the space $\Gamma$ of 
the representation $\gamma$ as a homothety with the coefficient 
\begin{equation}\label{koef}
(\lambda;\lambda+2\rho)=\|\lambda+\rho\|^2-\|\rho\|^2\;.
\end{equation}
From this formula it follows that
$$Tr(\gamma(b_{can}))=c_2(\gamma)\,Tr(\gamma(b_\gamma))=
(\lambda;\lambda+2\rho)\dim\Gamma\;.
$$
From the other hand, the definition of $b_\gamma$ implies that 
$$Tr(\gamma(b_\gamma))=\sum\limits_iTr(\gamma(e_i)\gamma(f_i))=
\sum\limits_i 1=\dim sl(r,\CC)=r^2-1\;.
$$
Therefore, we obtain:
$$
c_2(\gamma)=\frac{(\lambda;\lambda+2\rho)}{r^2-1}\dim\Gamma=
\frac{(\lambda;\lambda+2\rho)}{r^2-1}  
\prod_{\alpha\in R_+}\frac{(\lambda+\rho;\alpha)}{(\rho;\alpha)}\;.
$$
The last equality holds because of the formula for the dimension of 
the representation space with given highest weight [FH, 24.6]. This is the 
latter formula of the representation theory on which the calculation of 
$c_2(\gamma)$ is based. The Proposition is proved.

\vskip 10pt

In general case using (\ref{kasimir}) is 
hampered because the product by positive roots that equals  
$\dim\Gamma$ is cumbrous. However, for the exterior and symmetric powers 
of the canonical representation and for the adjoint representation $c_2$'s
are computed easily. 

\vskip 10pt

{\bf 1.3. Proposition.}

\vskip -10pt
\begin{equation}\label{c2la} 
c_2(\Lambda^k (can))={r-2\choose k-1}\ ,\quad
c_2(S^k (can))={r+k\choose k-1}\quad\mbox{\it and}\quad c_2(ad)=2r\;.
\end{equation}

{\sc Proof.} Consider the Cartan subalgebra
${\gtc h}\subset sl(r,\CC)$ consisting of diagonal matrices with zero 
trace. We shall define a matrix $H\in\gtc h$ by its diagonal elements
$(\lambda_1,\dots,\lambda_r)$.

The linear functions $\lambda_i-\lambda_j$, $i\ne j$, form the root system
$R\subset{\gtc h}^*$. The set
$S=\{\alpha_1,\dots,\alpha_{r-1}\}$, where
$\alpha_i=\lambda_i-\lambda_{i+1}$, is a basis in this system.
Matrices $H_i$ that have only two nonzero elements $\lambda_i=1$ and
$\lambda_{i+1}=-1$, $i=1,\dots,r-1$ form the dual basis in the dual system
$R^*\subset{\gtc h}$. Observe that $\alpha_i$ is identified with $H_i$ by
isomorphism ${\gtc h}^*\cong{\gtc h}$, which is defined by the Killing 
form
$(X;Y)=Tr(X\,Y)=B_{can}|_{\gtc h}(X;Y)$.
Herewith, the set $R_+$ of positive roots w.r.t. $S$, which  consists of 
the elements $\alpha=\lambda_i-\lambda_j$, $i<j$, is identified with the 
set of matrices $H_{\alpha}=H_i+\dots+H_j$, which also have two nonzero 
elements: $\lambda_i=1$ and $\lambda_j=-1$.

As it is known [FH], there is one-to-one correspondence between 
irreducible $sl(r,\CC)$-modules and their highest weights.
Herewith, if $\lambda$ is the highest weight of an irreducible 
representation $\gamma$, then $\gamma$ is finite dimensional iff
all values $\lambda(H_i)$, $i=1,\dots,r-1$,
are nonnegative integers.
This means that
$\lambda=k_1\omega_1+\cdots+k_{r-1}\omega_{r-1}$, where $k_i\in\ZZ$, 
$k_i\ge 0$ and $\omega_i$ are fundamental weights defined by the equalities
$\omega_i(H_j)=\delta_{ij}$. 

The fundamental weights have the form  $\omega_i =\lambda_1+\cdots +\lambda_i$,
$1\le i\le r-1$. We have:  $(\omega_i;\omega_k)=\frac{i(r-k)}{r}$ for
$i\le k$. 
It is easily seen that the half-sum of the positive roots is
$\rho =\sum_{i=1}^{r-1}\omega_i$. The highest weights of the 
representations $\Lambda^k (can)$, $S^k (can)$ and $ad$ are equal to
$\omega_k$, $k\omega_1$ and $\omega_1 + \omega_{r-1}$
respectively.  Calculating the homothety 
coefficients by formula (\ref{koef}), we get:  
$$(\omega_k\,;\,\omega_k+2\rho)=\frac{k(r-k)(r+1)}{r}\ ,\qquad
(k\omega_1\,;\,k\omega_1+2\rho)=\frac{k(r+k)(r-1)}{r}
$$
$$
(\omega_1 +\omega_{r-1}\,;\,\omega_1 +\omega_{r-1}+2\rho)=2\,r\,.
$$
By (\ref{kasimir}), multiplying these numbers respectively on
$\dim\Lambda^k (can)={r\choose k}$, $\dim S^k(can)={r+k-1\choose k}$, 
$\dim(ad)=r^2-1$ 
and then dividing by $r^2-1=\dim sl(r,\CC)$, we obtain the required formulae.
This completes the proof.

\vskip 10pt

Now we are ready to describe field configurations (i.~e. collections of
representations) such that the beta function coefficient is zero or negative.
As an answer we obtain the following statement.

\vskip 10pt

{\bf 1.4. Theorem.} {\parindent=0pt \it

1. For $su(2)$ one has $\beta(\{\gamma_i\})=0$ in the following cases:

a) $N_a=1$ and $\gamma_1=ad$.

b) all $\gamma_i=can$ and $N_a=4$;

The inequality $\beta(\{\gamma_i\})<0$ can be satisfied only in the case
b) while $N_a<4$.

2. For $su(3)$ one has $\beta(\{\gamma_i\})=0$ in the following cases:

a) $N_a=1$ and $\gamma_1=ad$.

b) all $\gamma_i$ are isomorphic either to $can$ or to $can^*\cong\Lambda^2 
(can)$, and $N_a=6$;

c) $N_a=2$, one of $\gamma_i$ is isomorphic either to $can\cong\Lambda^2 
(can^*)$ or $can^*\cong\Lambda^2 (can)$, and the other is either ${\rm 
S}^2(can)$ or ${\rm S}^2(can^*)$.

The inequality $\beta(\{\gamma_i\})<0$ can be satisfied either in the case
b) while $N_a<6$ or in a case of one representation isomorphic to
${\rm S}^2(can)$ or ${\rm S}^2(can^*)$.

3. For $su(r)$, $r\ge 2$, there is at least three series of cases in which
$\beta(\{\gamma_i\})=0$:

a) $N_a=1$ and $\gamma_1=ad$;

b) all $\gamma_i$ are isomorphic either to $can$ or to
$can^*\cong\Lambda^{r-1} (can)$, and $N_a=2r$. 

c) $N_a=2$, one of $\gamma_i$ is isomorphic either to $\Lambda^2 (can)$, 
or to $\Lambda^2 (can^*)$ and the other is either 
${\rm S}^2(can)$ or ${\rm S}^2(can^*)$.
}

\vskip 10pt

For $su(2)$ and $su(3)$ algebro-geometrical problem in some of these 
cases is beforehand divined and calculations are being performed using 
Vafa --- Witten's method (geometric vanishing situation).
For instance, the case 1a) is investigated in [VW], 1b) in [GL], 
algebro-geometrical problem for 2b) is considered in [Ku]. The case 2a) 
is geometrically analogous to 1a) (see below). At last, the most 
interesting case 2c) is a unique realization of the series 3c) on del Pezzo 
surfaces. The last section of this work is devoted to that.

\vskip 10pt

{\sc Proof of Theorem 1.4.}

The item 3 follows easily from Proposition 2. 
Suppose $r=2$. It is known that any non-trivial irreducible 
representation of $sl(2,\CC)$
is isomorphic to $S^k (can)$ for $k\ge1$. According to Proposition 2, we have:
$c_2(S^k (can))={k+2\choose k-1}={k+2\choose 3}$. This sequence of numbers is 
increasing. Besides, $ad\cong S^2 (can)$. This proves item 1.

To prove item 2, it suffices to estimate $c_2$ of $sl(3,\CC)$-representations. 
Suppose $\gamma$ is an irreducible finite dimensional 
$sl(3,\CC)$-representation, $P$ 
is the set of its weights, and $\Gamma$ is the representation space. For any $\pi\in P$ denote by
$\Gamma^{\pi}$ the corresponding weight subspace. 
The number
$m_\pi=\dim \Gamma^{\pi}$ is called the multiplicity of $\pi$. It is known 
that $P$ is invariant under the action of the Weyl group $W$, which is generated 
by root reflections, and for any $\pi\in P$ and $w\in W$ the 
multiplicities of weights $\pi$ and $w(\pi)$ coincide. Besides,
$$\Gamma=\bigoplus_{\pi\in P}\Gamma^{\pi}\;.$$
 By definition, $c_2(can)=1$ and 
it is evident that for the canonical representation $Tr(H_1^2)=2$. Let
$\lambda$ be the highest weight of $\gamma$, then by (\ref{c2viabxy}) 
we have:
$$c_2(\gamma)=\frac{1}{2}\,Tr(\gamma(H_1)^2)=
  \frac{1}{2}\sum_{\pi\in P}\pi(H_1)^2m_\pi\ge
  \frac{1}{2}\sum_{\pi\in\lambda^W} \pi(H_1)^2\,.
$$
Here $\lambda^W$ is the orbit of the highest weight under the action of the
Weyl group and $m_\lambda$ always equals 1.
If $\lambda=k_1\omega_1+k_2\omega_2$,  then the last half-sum can 
be easily calculated. It is equal to $2(k_1^2+k_2^2)$ for $k_1,k_2>0$ and 
$k^2$ if one of $k_1,k_2$ equals 0 and the other is $k$.
Suppose $\gamma$ is contained in a collection of representations such that
the beta function coefficient is non-positive; then 
$c_2(ad)=6\ge c_2(\gamma) $, whence
in the former case $k_1^2+k_2^2\le 3$ i.e. $k_1=k_2=1$ and $\gamma=ad$ and 
in the latter case there are four possibilities:
$$\renewcommand{\arraystretch}{1.5} 
  \begin{array}{|c|c|c|c|c|}
    \hline
  (k_1,k_2)& (1,0) &    (2,0)      & (0,1) & (0,2)\\ \hline
   \gamma    &  can  &{\rm S}^2(can) & can^* & {\rm S}^2(can^*)\\ \hline
 c_2(\gamma) &   1   &       5       &    1  &        5        \\  \hline
 \end{array}
$$
This concludes the proof.

\vspace{2ex}

{\large\bf 2. From field equations to algebro-geometrical problems.}

\nopagebreak[3]

\vspace{2ex}

Theorem 1 shows that we have three series of field configurations:
a), b) and c), 
which take place for all $su(r)$. We consider only the first two of them 
in this section. 

{\bf 2.1. }For series a) our fields consist of multiplets in the adjoint
representation. We write field equations for a special type of metrics, 
namely K\"ahler metrics, because our calculations will be related only to 
the case of algebraic surfaces, where necessary for us constants will have 
algebro-geometrical meaning and so can be computed within the framework 
of algebraic geometry. For simplicity we shall consider simply connected 
case. So, let $S$ be simply connected K\"ahler surface and $K$ be its 
canonical class. Configuration space is the direct product
$$
{\cal A}(E)\times\Gamma(ad E\otimes K)
$$
where ${\cal A}(E)$ is affine space of $SU(r)$-connections in a
complex vector bundle $E$ on $S$ of rank $r$ with $c_1=0$ and $c_2=k$. 
$\Gamma(ad E\otimes K)$ is the space of Higgs --- Hitchin fields.
Then for a point
$(a,\phi)\in {\cal A}(E)\times\Gamma(ad E\otimes K)$
the field equations have the form
\begin{equation}\label{pe1}
\arraycolsep=0.1em
\begin{array}{rcl}
\displaystyle{F_a^{2,0}=F_a^{0,2}}&=&0 \medskip \\
\displaystyle{\overline{\partial}_a (\phi)}&=&0 \medskip \\
\displaystyle{\omega\wedge F_a+[\phi,\overline{\phi}] }&=&0
\end{array}
\end{equation}
where $F_a^{i,j}$ are the Hodge decomposition components for the 
curvature tensor 
$F_a$ of the connection $a$ and $\omega$ is the K\"ahler form of our metric.
Vafa and Witten in the work [VW] wrote generalizations of the equations above
for a generic metric.
Using linearization of these equations it is easy to show that the space 
of solutions is oriented and finite. Hence, for any $k=c_2(E)$
we obtain a number $a_k$ of these solutions with provision for the 
orientation sign.
Our partition function in this case has the form
$$
Z_{\mbox{\footnotesize (a)}}=\sum_{k=0}^{\infty}a_k q^k
$$
where as usual $q=e^{2\pi i k\tau}$ (see details in [VW]).

Geometric meaning of equations (\ref{pe1}) is as follows.
The first equation is equivalent to statement that the curvature form has 
type (1,1), this is equivalent to defining a holomorphic structure on 
$E$. The second equation is equivalent to holomorphicity of the 
map $\phi:E\longrightarrow E\otimes K$ and the third equation has a sense 
of zero level of
the moment map
for the gauge group action, i.e. it means some stability conditions 
for the holomorphic pair 
$(E\,,\,\phi)$.

Now suppose that for our K\"ahler structure on $S$ this stability condition implies 
vanishing of the holomorphic Higgs --- Hitchin field $\phi$. It will be so, 
for example, if the canonical class $K\le 0$, i.e. if $S$ is  
a del Pezzo surface or K3 surface. In this case our K\"ahler metric is 
non-general because the system
(\ref{pe1}) comes to anti-self-duality (ASD) equations on the connection $a$.
We use here the following Donaldson's result [DK]: the set of gauge orbits
of irreducible
ASD connections on the bundle $E$ is in one-to-one correspondence with the set 
of stable holomorphic structures on $E$, where the stability is by 
Mumford --- Takemoto w.r.t. $[\omega]$.
Therefore, the moduli 
space of solutions up to gauge of the field equations (\ref{pe1})
is the whole space ${\cal M}_k(r)$ of $SU(r)$-instantons on
$S$ with the charge $k$. For example, for $r=2$ this variety has dimension
$4k-3$ for del 
Pezzo surfaces and $4k-6$ for K3 surfaces. To regularize the problem we
are to apply deformation to normal cone as described in chapter III of [PT].
In this case the obstruction bundle on ${\cal M}_k(r)$ has as a fiber the 
cokernel of the map $\overline{\partial}_a $, i.e. over a point 
$E\in {\cal M}_k(r)$ the fiber of the 
obstruction bundle is $H^1(adE\otimes K)$ by  
Dolbeault   
isomorphism. Thus, the obstruction bundle for our problem is 
$\Omega {\cal M}_k(r)$
by Kodaira --- Spenser
isomorphism. According to construction expounded in chapter III of [PT] the number 
of  solutions of our problem is
\begin{equation}\label{e4}
a_k=c_{top}(\Omega {\cal M}_k(r))
\end{equation}
i.e. the top Chern class of the cotangent bundle.

\vspace{1ex}

{\bf Remark.} Variety ${\cal M}_k(r)$ can be non compact (for instance, this takes 
place for $r=2$) and we should consider it as a variety "with ends". We are to 
estimate which number of solutions approaches infinity on these "ends" from general
number of solutions  defining by topological Euler characteristics of 
Fried --- Uhlenbeck
compactification of ${\cal M}_k(r)$ considering as orbifold (see [DK]).
(Unfortunately, this analysis is absent in basic work [VW], and the reader
is to make it himself using as a model analogous problem of expressing 
spin-polynomials in terms of Donaldson polynomials. A first
step of solving this
problem is described in [W-L]).

\vspace{1ex}

Thus, (after realization of the program expounded in the remark)  the partition 
function  for a) series takes the form
\begin{equation}\label{e5}
Z_{\mbox{\footnotesize (a)}}=\sum_{k=0}^{\infty}\chi(\overline{{\cal M}_k(r)}) q^k\,,
\end{equation}
where  $\chi$ is the topological Euler characteristics of
$\overline{{\cal M}_k(r)}$ as an orbifold.

Testing this sum for modularity performed in  [VW] (for K3 surfaces and $\CC\PP^2$)
is based on a remarkable work [Kl] of Klyachko who has expressed coefficients 
of the
series (\ref{e5}) in terms of arithmetical 
Hurwitz 
function.

{\bf 2.2. } We now describe a physical setup for the b) series. 

Let $M,g$ be smooth compact simply-connected Riemannian manifold. Consider the 
Levi --- Civit{\'a } connection on the tangent bundle. A choice of 
$Spin^{\Bbb C}$-structure gives us a pair of complex Hermitian bundles 
$W^{\pm} $ of rank 2 with a class $c = c_1(\det W^{\pm})$.

The Levi --- Civit{\'a } connection induces an $SO(3)$-connections on  
$W^{\pm}$.
Consider a Hermitian bundle $E$ of rank $r$ with unitary
connection $a$ defined up to homothety.
Besides, consider a $U(1)$-connection $b$ on line bundle $\det(E\otimes W^+)$.

Configuration space for our field equations is the set of 
collections 
$$
(a,b,\phi_1,\dots,\phi_{2r})\,,
$$
where  $\phi_j$ are sections of the tensor product
$E\otimes W^+$.

Recall that Dirac operator $D_{a,b}$ is a differential operator
defined by the following way. Covariant derivative determines 
a map
\begin{equation}\label{cdif}
d_{a,b}: \Gamma(E\otimes W^+) \longrightarrow \Gamma(E\otimes W^+
\otimes T^*M),
\end{equation}
and the last tensor product is contracted as
$$\Gamma(E\otimes W^+ \otimes  T^*M) = \Gamma(E\otimes W^+ \otimes
\hom(W^-,  W^+)^*)   =   \Gamma(E\otimes W^-)$$
by means of Clifford multiplication. 
Superposition of covariant derivation and contraction gives us the twisted 
Dirac operator 
\begin{equation}\label{dirac} 
D_{a,b}: \Gamma(E\otimes W^+) \longrightarrow \Gamma(E\otimes W^-)\;.  
\end{equation}
Our field equations take the form:  
\begin{equation}\label{pee} 
    \begin{array}{l} 
\displaystyle{D_{a,b}(\phi_j) = 0} \medskip \\
\displaystyle{F_a^+     = i\sum_j  (\phi_j\otimes
\overline{\phi_j})_{0}}\ ,
\end{array}
\end{equation}
where the index $_{0}$ means taking traceless part in
$\End\,(W^+)$. Here we use the natural identification $ad\,W^+\cong \Omega^+$ 
(see [T]).

\vspace{1ex}

{\bf Remark.} One can see that these equations have the form of zero 
level for  
moment map
of hyper-K\"ahler reduction w.r.t. the gauge group.
In particular, the first equation contains an interpretation of zero 
level of  
moment map
for holomorphic symplectic structure.

\vspace{1ex}

Now suppose that we have chosen a bundle $E$ such that moduli space 
${\cal M}_c(E)$ of orbits 
of  solutions of (\ref{pee}) is finite. Let $\#({\cal M}^g_c(E))$ be the
number of solutions weighted with the orientation sign. Since 
$\#({\cal M}^g_c(E))$ depends only on topological type of the bundle
$E$, we can write
$$\#({\cal M}^g_c(E))=\#({\cal M}^g_c(r,c_1,k))\ ,$$
where as usual $k=c_2(E)$.

Our partition function has the form
$$
Z_{\mbox{\footnotesize (b)}}=
\sum_{k=0}^{\infty}\#({\cal M}^g_c(r,c_1,k)) q^k\ ,
$$
where again $q=e^{2\pi i\tau}$ and the metric 
$g$ is sufficiently general (for non-general metrics the moduli space can 
have parameters).

\vspace{1ex}

{\bf Remark.} It is easy to see that our equations define a multiplet of
non-abelian mo\-no\-po\-les which are generalizations of abelian Seiberg ---
Witten monopole (see, for example, [T]).

\vspace{1ex}

Recall that we started with a K\"ahler metrics in describing a) series i.e.
we reduced our problem to pure algebro-geometrical problem. To do this 
for b) series we have to act more accurately. Here we illustrate the general 
situation by the case 
$r=2$ and $M=\CC\PP^2$.

Consider the standard Fubini --- Study metric $g_{FS}$.
This metric has positive scalar curvature. The Weitzenb\"ok
formula shows that in this case the multiplet of twisted harmonic spinors
$(\phi_1,\dots,\phi_{4})$ vanishes, 
and the field equations come to ASD equations. The moduli space of
solutions ${\cal M}^{g_{FS}}_c(2,c_1,k)$ coincides with the moduli space of
instantons (holomorphic stable bundles) of this topological type. 
Arithmetic calculations of the next 
section show that the unique possibilities for $c_1$ are $-1$ or $-5$ and
$$
\dim {\cal M}^{g_{FS}}_c(2,-1,k)=4(k-1)\,,
\quad\dim {\cal M}^{g_{FS}}_c(2,-5,k)=4(k-7)\;.
$$

To regularize the problem we again have to apply the deformation to
the normal cone as it is described in [PT, chapter III].
In this case the fiber of the obstruction bundle over 
a point in ${\cal M}^{g_{FS}}_c(2,c_1,k)$ 
is the  direct sum of four cokernels of twisted Dirac operator
(\ref{dirac}). As usual, we interpret the an ASD connection as a holomorphic
stable structure on $E$ and Dirac operator as 
$\overline{\partial}_a \oplus {\overline{\partial}_a }^*$, where 
$Spin^{\Bbb C}$-structure equals $-K_S$ 
(see, for example, [DK]).
So, the cokernel of Dirac operator is the coherent cohomology
space  $H^1(E)$.

Thus, our moduli spaces ${\cal M}^{g_{FS}}_c(2,c_1,k)$ are the moduli spaces of 
stable bundles ${\cal M}(2,c_1,k)$
with this topological type, where $c_1\in\{-1,-5\}$. 
Hence, the fiber of the obstruction bundle over a point
$E\in {\cal M}^{g_{FS}}_c(2,c_1,k)$ is $H^1(E)^{\oplus 4}$.
According to the construction given in chapter III of [PT] the number of 
 solutions of our problem is
\begin{equation}\label{e4b}
\#({\cal M}^{g_{FS}}_c(2,c_1,k))=c_{top}(({\cal H}_k)^{\oplus 4})\;,
\end{equation}
where ${\cal H}_k$ is a bundle over the moduli space with the fiber $H^1(E)$
and specially chosen normalization (see [GL]).

In fact the above reasoning hold when $S$ is a del Pezzo surface,   
$Spin^{\Bbb C}$-structure is $-K_S$ and $g$ is the Hodge metric. In this situation the space ${\cal M}^g_{-K_S}(r,c_1,k)$ of solutions up to gauge
of the field equations is also interpreted as 
the moduli space ${\cal M}(r,c_1,k)$ of 
stable bundles with the given topological type. The fiber of the 
obstruction bundle over a point $E$ is
the direct sum of $2r$ cokernels of the Dirac operator i.e.
$H^1(E)^{\oplus 2r}$. Hence, to have a finite number of solutions it is necessary 
to have the equality
$$
\dim{\cal M}(r,c_1,k)=2r\dim H^1(E)\;.
$$
If this condition holds, then the partition function is

\begin{equation}\label{zb}
Z_{\mbox{\footnotesize (b)}}=
\sum_{k=0}^{\infty}\#({\cal M}^g_{-K_S}(r,c_1,k)) q^k =
\sum_{k=0}^{\infty}c_{top}(({\cal H}_k)^{\oplus 2r}) q^k  
\end{equation}
where ${\cal H}_k$ is a bundle over ${\cal M}(r,c_1,k)$ with the 
fiber $H^1(E)$.

Observe that the moduli space ${\cal M}(r,c_1,k)$ is non-compact as a 
rule, so $c_{top}$'s in the last formula depend on cohomology theory
which we are working with. A natural way to give a meaning to this
quantities is to consider the Gieseker compactification
$\overline{\cal M}(r,c_1,k)$ which is the coarse moduli space of 
semistable sheaves. We discuss this in the next section.

\vspace{2ex}

{\large\bf 3. Exceptional bundles in realization of b) series.}

\nopagebreak[3]

\vspace{2ex}

{\bf 3.1.} Let $S$ be a del Pezzo surface with the canonical class $K$, 
${\cal M}(r,c_1,c_2)$ be the moduli space of
stable bundles over $S$ with the given topological invariants, 
$\overline{\cal M}(r,c_1,c_2)$ be the Gieseker compactification of
${\cal M}(r,c_1,c_2)$. The stability notion is by Mumford --- Takemoto w.r.t.
$(-K)$. This means that a sheaf $E$ is stable if it is torsion free and
for any subsheaf $F$ of $E$ with $r(F)<r(E)$ 
$$
\mu(F)<\mu(E)\,,\quad\mbox{where}\quad\mu=\frac{c_1\cdot (-K)}{r}\quad 
\mbox{is the slope.}
$$
In this section we regard ${\cal M}(r,c_1,c_2)$ as the space 
 of solutions of (\ref{pee}) up to gauge. We
consider the partition function (\ref{zb}) under the 
following assumptions: 
\begin{equation}\label{dim}
\hbox to 0.9\textwidth{%
\mbox{\it dimension condition}\hfil%
$\dim{\cal M}(r,c_1,k)=2r\dim H^1(E)$\hfil}
\end{equation}
\begin{equation}\label{mub}
\hbox to 0.9\textwidth{%
\mbox{\it slope inequality}\hfil%
$-K^2<\mu(E)<0$\hfil}
\end{equation}

The dimension condition is necessary for the coefficients of the partition
function (\ref{zb}) to be really a numbers. From the slope inequality,
the stability properties and Serre duality it follows that
$H^0(E)=H^2(E)=0$. Therefore, we have:
\begin{equation}\label{h1} 
\dim H^1(E)=-\chi(E)=k-\frac{1}{2}c_1\cdot (c_1-K)-r\;.
\end{equation}
Further, the space $\ext^1(E,E)$ coincides with formal tangent space
to ${\cal M}(r,c_1,k)$ at the point $E$. Using stability of $E$, we have
$\hom(E,E)\cong\CC$. In addition, $(-K)$ is ample, so by Serre duality 
$\dim\ext^2(E,E)<\dim\hom(E,E)$ (see [KO]). This implies that
$\ext^2(E,E)=0$, and we obtain
\begin{equation}\label{dm}
\dim{\cal M}(r,c_1,k)=1-\chi(E,E)=1-\chi(E^*\otimes E)=2rk-(r-1)c_1^2-r^2+1
\end{equation}

\vspace{1ex}

{\bf Remark.} If the numbers $r$ and $c_1\cdot K$ are coprime, then any
semistable sheaf in $\overline{\cal M}(r,c_1,c_2)$ is stable. This implies 
that the moduli space $\overline{\cal M}(r,c_1,c_2)$ is smooth, since
$\dim\ext^1(E,E)$ is independent of $E$.

\vspace{1ex}

By the equalities (\ref{h1}) and (\ref{dm}), the dimension
condition (\ref{dim}) takes the form:
\begin{equation}\label{equb}
c_1^2-r\,c_1\cdot K+r^2+1=0
\end{equation}
The most interesting  solutions of this equation are on the complex projective plane.

\vspace{1em}

{\bf 3.2. Theorem.} {\it A stable sheaf $E$ on $\PP^2$ with invariants satisfying
(\ref{mub}) obeys the equality (\ref{dim}) iff the triple of numbers
$(r,-c_1,1)$ satisfies the Markov equation
\begin{equation}\label{Markov} x^2+y^2+z^2=3\,xyz\ . \end{equation} }

The proof follows immediately from the fact that $c_1\cdot K=-3\,c_1$ for $\PP^2$.

\vspace{1em}

The solution set of (\ref{Markov}) is well-known [GR]. Namely, any solution
can be obtained from $(1,1,1)$ by a finite number of mutations. A mutation
is changing one 
variable to another root of the quadratic equation while 
the rest variables are fixed. So, for $z=1$ we have a sequence of mutations
that gives the set of all solutions $(r,c_1)$ to (\ref{equb}). Here are
several first solutions such that $r>-c_1$:
\begin{equation}\label{seq}
(2,-1)\,,\ (5,-2)\,,\ (13,-5)\,,\ (34,-13)\,,\ (89,-34)\,,\
(233,-89)\,,\ \dots 
\end{equation}
The pair after $(r,c_1)$ in this sequence is $(3r+c_1,-r)$. 

It is known [GR] that for any exceptional collection of three bundles on
$\PP^2$ the ranks form a  solution of Markov equation. Moreover, there
are natural mutations of exceptional collections, which induce the
above-mentioned mutations of rank collections. Let us enumerate the pairs 
in (\ref{seq}); then there is a unique sequence $E_n$ of exceptional bundles 
such that $(r(E_n),c_1(E_n))$ is exactly the $n$-th term of (\ref{seq}).
For example, $E_1\cong\Omega_{\PP^2}(1)\cong\mbox{T}\PP^2(-2)$. 
Moreover, for any $n$ collection $(E_n,E_{n+1},{\cal O}_{\PP^2})$ 
is exceptional and 
the mutations preserving ${\cal O}_{\PP^2}$ are given by the exact triples 
\begin{equation}\label{mut}
0\longrightarrow E_{n-1}\longrightarrow \hom (E_n,E_{n+1})\otimes E_n
  \longrightarrow  E_{n+1} \longrightarrow 0\,,
\end{equation}
where $\hom (E_n,E_{n+1})$ is canonically identified with 
$\hom (E_{n-1},E_{n})^*$. In fact, $E_n$'s are defined also for $n\le 0$ by
the  triple above.
For example, $E_0={\cal O}_{\PP^2}(-1)$,  $E_{-1}={\cal O}_{\PP^2}(-2)$,
and for $n=0$ the exact triple (\ref{mut}) is the Euler exact sequence
twisted by $(-2)$. The bundles 
$E_n=\overline{\cal M}(r_n,-r_{n-1},c_2(E_n))$ 
are the initial points for series of moduli spaces  
$\overline{\cal M}(r_n,-r_{n-1},k)$, where $k\ge c_2(E_n)$. 
Any of these spaces is smooth because $r_n$ and $3\,r_{n-1}$ are coprime.
Thus, we have a diagram:
$$\renewcommand{\arraystretch}{2}
\begin{array}{ccccl}
\vdots &\vdots &\vdots  \\
{\cal O}_{\PP^2}(-2)&
       Hilb^1(\PP^2)&
            Hilb^2(\PP^2)&\dots &Z_{-1}\mbox{ is a polynomial of degree }\le4\\

{\cal O}_{\PP^2}(-1)&
      \stackrel{\raisebox{1.5ex}{$|\,|$}}{\PP^2}&
            Hilb^2(\PP^2)&\dots\ &Z_0=1+q\\

\Omega_{\PP^2}(1)&
      \overline{\cal M}(2,-1,2)&
            \overline{\cal M}(2,-1,3)&\dots\ &
Z_1=q+13q^4+729q^5+85025q^6+
\\
E_2        &\overline{\cal M}(5,-2,5)
              &\overline{\cal M}(5,-2,6)&\dots\ &
            \qquad\qquad\qquad\qquad\raisebox{2ex}{$+15650066q^7+\cdots$}\\
E_3        &\overline{\cal M}(13,-5,19)
              &\overline{\cal M}(13,-5,20)&\dots\\
\vdots &\vdots &\vdots &\ddots \\
\end{array}
$$
Each row of this diagram produces a partition function $Z_n$. 
All what is known about these functions for the moment is written
in the right hand side.
On the other hand,
the terms of the first column are connected by (\ref{mut}). 

\vspace{1ex}

{\bf Question.}
Exist there a transformation between the whole rows of the diagram
that generalizes the mutations (\ref{mut})?

\vspace{1ex}

{\bf Conjecture.} If such a transformation exists, then one can calculate
all the partition functions knowing only two of them.

\vspace{1ex}

We would like to observe that the calculation of the first coefficients 
of $Z_1$ is contained in the paper [GL]. It is important that the 
moduli spaces $\overline{\cal M}(2,-1,k)$ are really fine i.e. there
exist universal families which are needed to construct the bundles 
${\cal H}_k$.  
The first part of the paper [GL] is description of $\overline{\cal M}(2,-1,k)$
as fine moduli spaces.  A similar description
of all $\overline{\cal M}(r_n,-r_{n-1},k)$ in the diagram above 
and some moduli spaces of sheaves over $\PP^1\times\PP^1$ is given
in [Ku2]. 
The second part of [GL] explains the calculations of the coefficients via
Bott's formula performed by MAPLE V on computer.

{\bf 3.3.} In fact, the situation described above for $\PP^2$  admits a
generalization for del Pezzo surfaces. The bundles $E_n$ described in the 
previous subsection are exactly all exceptional bundles that belong to
the subcategory
$$
{\cal O}^\perp=\{F:\ \ext^i({\cal O},F)=0\}
$$
of the category of coherent sheaves. We see that for $\PP^2$ each 
exceptional bundle from ${\cal O}^\perp$ (i.e. with zero cohomology
groups) produces a partition function. It is not hard to check that 
the slope inequality (\ref{mub}) holds for these bundles automatically. 
That is not so for other del Pezzo surfaces.
For example, the bundles ${\cal O}(-1, m)$ over $\PP^1\times\PP^1$ belong to
${\cal O}^\perp$ but have the slope $2\,(m-1)$.

\vspace{1ex}

{\bf Proposition.} {\it If an exceptional bundle $E$ over del Pezzo surface 
satisfies the slope inequality (\ref{mub}), then} 
$$
E\in{\cal O}^\perp\ \Longleftrightarrow\ 
\mbox{\it the moduli spaces }\overline{\cal M}(r(E),c_1(E),k)\ 
\mbox{\it obey the dimension condition\/.} 
$$
{\sc Proof.} $(\Rightarrow)$ Obviously, the dimension condition is
independent of $k$ and holds for $E=\overline{\cal M}(r(E),c_1(E),c_2(E))$.

\noindent $(\Leftarrow)$ 
It is known [KO] that exceptional bundles are stable
by Mumford --- Takemoto w.r.t. $(-K)$. Therefore, as in {\bf 3.1}
(\ref{mub}) implies that $H^0(E)=H^2(E)=0$. By the dimension condition
applied to $E=\overline{\cal M}(r(E),c_1(E),c_2(E))$, we have:
$H^1(E)=0$. This concludes the proof.

\vspace{1ex}

{\bf Corollary.} {\it Each exceptional bundle from ${\cal O}^\perp$
that satisfies the slope inequality gives a partition function for
mathematical testing of $S$-duality conjecture.}

\vspace{1ex}

In particular, any system of bundles generated by an exceptional pair
produces a diagram of moduli spaces which is similar to the diagram
for $\PP^2$. The question and the conjecture are the same.

\vspace{2ex}

{\large\bf 4. An algebro-geometrical model for c) series.}

\nopagebreak[3]

\vspace{2ex}

Physical setup for the c) series in unknown for the moment.
In this section we consider pure mathematical problem which is similar
to one in the previous section. 

Suppose that $E$ is stable bundle with $r\ge 3$ 
over a del Pezzo surface and
\begin{equation}\label{muv}
-\frac{1}{2}K^2<\mu(E)<0\;.
\end{equation}
The spaces $H^0(S^2E)$ and  $H^0(\Lambda^2E)$ are canonically identified
with the spaces of symmetric and skew-symmetric maps of bundles
$E^*\longrightarrow E$. From stability properties and the inequality
$\mu(E)<0<\mu(E^*)$ it follows that $\hom(E^*,E)=0$. Therefore, we have: 
$H^0(S^2E)=H^0(\Lambda^2E)=0$.

By Serre duality, $H^2(S^2E)^*=H^0(S^2E^*\otimes K)$. The last space 
identifies with the space of maps
$\varphi:E\longrightarrow E^*\otimes K$ such that
$\varphi^*(K)=\varphi$. Since $0<\mu(E^*)<\frac{1}{2}K^2$, we have 
$\mu(E^*\otimes K)<-\frac{1}{2}K^2<\mu(E)$. This implies that
$\hom(E,E^*\otimes K)=0$ and so, $H^2(S^2E)=0$. By similar reasoning, we
have $H^2(\Lambda^2E)=0$. 
Thus, the dimensions of $H^1$ spaces for symmetric and exterior square of
the bundle $E$ are expressed via $(r,c_1,k)$ by the Riemann --- Roch formula:  
$$\dim H^1(S^2E) = -\chi(S^2E)= 
(r+2)k-\frac{r+3}{2}c_1^2+\frac{r+1}{2}c_1\cdot K-\frac{r(r+1)}{2}\;, $$ 
$$\dim H^1(\Lambda^2E) = -\chi(\Lambda^2E) =
(r-2)k-\frac{r-1}{2}c_1^2+\frac{r-1}{2}c_1\cdot K-\frac{r(r-1)}{2}\;. $$ 

Taking into account (\ref{dm}) we conclude that the dimension condition
\begin{equation}\label{dmv} 
\dim {\cal M}(r,c_1,k)=\dim 
H^1(S^2E)+\dim H^1(\Lambda^2E) 
\end{equation} 
is equivalent to the equation on $r$ and $c_1$:  
\begin{equation}\label{eqv} 
2\,c_1^2-r\,c_1\cdot K+1=0 \;.
\end{equation}
It is easy to see that for $\PP^2$ this equation has no solutions. Indeed, 
in this case we have $c_1\cdot K=-3c_1$ and discriminant of  (\ref{eqv})
w.r.t. $c_1$ equals $9r^2-8$. This can not be a square of an integer
for $r\ge 3$ because the difference between neighbouring squares is
$9r^2-(3r-1)^2=6r-1>16$.

For $\PP^1\times\PP^1$ there are no solutions also because
$c_1\cdot K$ is even.

Now consider the plane with $m$ blowed points in general position, 
where $1\le m\le 8$. Denote this surface by $X_m$.
We have:
$$\Pic X_m=\ZZ\ell_0\oplus\ZZ\ell_1\oplus\ZZ\ell_2\oplus\cdots\oplus
\ZZ\ell_m\;,
$$
where $\ell_0$ is preimage of a line in $\PP^2$ and $\ell_i$, $1\le i\le m$, 
are exceptional curves.
The canonical class of $X_m$ is
$$K=-3\ell_0+\sum\limits_{i=1}^m\ell_i\ ,
$$
and the intersection form is defined by relations:
$$\ell_0^2=1\,,\ \ell_i^2=-1\,,\ 1\le i\le m\,;\quad 
\ell_i\cdot\ell_j=0\,,\ i\ne j\;.
$$
Consider $c_1=a\ell_0+\sum\limits_{i=1}^m b_i\ell_i$. We have
$c_1\cdot K=-3a-\sum\limits_{i=1}^m b_i$;
$c_1^2     = a^2-\sum\limits_{i=1}^m b_i^2$.
Suppose $c_1\cdot K=A$; then from (\ref{eqv}) we obtain:
$c_1^2=\frac{1}{2}(rA-1)$. Together with (\ref{muv}) this gives the
system of conditions:
\begin{equation}\label{sys}
\left\{\begin{array}{l} 
   \sum\limits_{i=1}^m b_i+3a+A=0\smallskip \\                          
    \sum\limits_{i=1}^m b_i^2 =a^2-\frac{1}{2}(rA-1)\smallskip \\
         0< A < \displaystyle{\frac{r}{2}(9-m)}\ .
\end{array}\right.
\end{equation}
The first equation defines a hyperplane in the Euclidean space with coordinates
$(b_1,\dots,b_m)$ and the second one defines a sphere with center in the
origin. Now we show that in all cases except only one the distance  
$d=\frac{|3a+A|}{\sqrt{m}}$ between the origin end the hyperplane is greater
then the radius $R$ of the sphere. Consider
$$m(d^2-R^2)=(9-m)a^2+6a+A^2+\frac{m}{2}(rA-1)\,.
$$
This is a quadratic polynomial w.r.t. $a$ with discriminant
$$D=\frac{m}{2}(rA-1)\left(m-9+\frac{2A^2}{rA-1}\right)\ .
$$
For $0<m<9-\frac{2A^2}{rA-1}$ we have $D<0$, so $d^2-R^2>0$, and
there are no solutions of the system (\ref{sys}). Hence, if there are
solutions, then the number of blowed points and $c_1\cdot K$ are restricted 
with the inequality:
$m\ge 9-\frac{2A^2}{rA-1}$. 
From the other hand, the third inequality of (\ref{sys}) implies that
$m<9-\frac{2A}{r}$, i.e. $m$ belongs to the interval
$$I=\left[9-\frac{2A^2}{rA-1}\;;\ 9-\frac{2A}{r}\right)\ .$$ 
Since $A$ and $r$ are integers, we have that the right bounds of $I$ are
either integers or 
distant
from integers by the distance greater than or equal to $\frac{1}{r}$. At the
same time the length of  $I$ is
$\frac{2A^2}{rA-1}-\frac{2A}{r}=
\frac{2}{r(r-\frac{1}{A})}\le\frac{2}{r(r-1)}$, that is less then or equal to
$\frac{2}{3r}$ for $r\ge 4$. Therefore, for $r\ge 4$ the interval $I$ 
does not contain an integer and the system (\ref{sys}) has no solution.
For  $r=3$ after a simple computation we see that there is a unique solution:
$m=8$, $A=1$,
$a=-3$, $b_1=\cdots=b_8=1$. Thus, we have proved the following

\vskip 10pt

{\bf Proposition. }{\it A stable sheaf $E$ over $X_m$ satisfies
to the restriction (\ref{muv}) and the dimension condition (\ref{dmv}) iff
$m=8$, $r(E)=3$ and $c_1(E)=K$.}

\vskip 10pt

Observe that under conditions of the Proposition
$c_1(E)\cdot K$ and $r(E)$ are coprime, so the moduli space of stable bundles
in question  has smooth Gieseker compactification. Therefore, 
we have an algebro-geometrical problem ("for testing $S$-duality
conjecture"). Namely, let
${\cal S}_k$ and ${\cal L}_k$ be bundles over ${\cal M}(3,K,k)$ with
fibers $H^1(S^2E)$ and $H^1(\Lambda^2E)$. The problem is to test the 
function
$$Z=\sum_kc_{top}({\cal S}_k\oplus {\cal L}_k)\,q^k\ .
$$
for modular behaviour.

\vspace{2em} 

\centerline{\large\bf References.}

\nopagebreak[3]

\vspace{1em} \parindent=0pt

[APS] Ph.C.Argyres, M.R.Plesser, N.Seiberg.  The Moduli Space of Vacua of 
$N=2$ SUSY QCD and Duality in $N=2$ SUSY QCD //Nucl. Phys. B {\bf  471}
(1996). P.159-194. (available also from hep-th/9603042)

\vspace{1em}

[GR] A. L. Gorodentsev, A.~N.~Rudakov. Exceptional Vector Bundles on
    Projective Spaces// Duke Math.~J., 54 (1987) 115--130.

\vspace{1em}

[GL] A.L.Gorodentsev, M.I.Leenson.   How to calculate the correlation 
function in twisted SYM $N=2$, $N_f=4$ QFT on projective plane// Preprint 
MPI 96-49

\vspace{1em}

[DK] S. Donaldson and P. Kronheimer. The Geometry of Four-Manifolds//
Clarendon Press,
Oxford,
1990.

\vspace{1em}

[DW] R.Donagi, E.Witten.  Supersymmetric Yang-Mills theory and integrable
systems// Nucl. Phys. B {\bf  460}
(1996). P. 299-334.(available also from  hep-th/9510101)

\vspace{1em}

[FH] W. Fulton, J. Harris.  Representation Theory. A First Course//   
Springer-Verlag (1990).

\vspace{1em}

[Kl] A.~A.~Klyachko. Moduli of vector bundles and numbers of classes//
Translation: Fun\-ctio\-nal-Anal-Appl. {\bf 25} (1991), no.~1. P.~67--69. 

\vspace{1em}

[KO] S.~A.~Kuleshov and D.~O.~Orlov. Exceptional sheaves on del Pezzo
surfaces// Russian Acad. Sci. Izv. Math. Vol.~{\bf 44} (1995) No~1. P.~479--513.

\vspace{1em}

[Ku1] S.A.Kuleshov.   Moduli Space of Bundles on a Quadric//
    Warwick Preprints: 18/1996

\vspace{1em}

[Ku2] S.A.Kuleshov.   Moduli Spaces of sheaves necessary for
  testing $S$-duality conjecture// Preprint MPI 97-32.

\vspace{1em}

[MO] C.~Montonen, D.~Olive, Phys. Lett 72B (1977), 117-132.

\vspace{1em}

[PT] V.~Ya.~Pidstrigach and A.~N.~Tyurin. Invariants of the smooth 
structure of an algebraic surface arising from the Dirac operator//
Russian Acad. Sci. Izv. Math. Vol. {\bf 40} (1993) No~2. P.~267--351.

\vspace{1em}

[T] N.~A.~Tyurin. Necessary and sufficient condition for a deformation
of a B-monopole into an instanton// Izvestiya: Mathematics, Vol.~{\bf 60}
(1996), no.~1. P.~217--231. 

\vspace{1em}

[VW] G.~Vafa, E.~Witten. A strong coupling test of S-duality//
    Nucl. Phys. B {\bf  431}
(1995). P.3--77 (available also from hep-th/9408074).

\vspace{1em}

[W] E.Witten. Monopoles and Four-Manifolds// Math. Res. Letts. 1 (1994),
    769--796.

\vspace{1em}

[W-L] W.-M.~R.~Leung. On {\it Spin}$^\CC$-Invariants of Four-Manifolds//
    Thesis PhD, Magdalene College, Oxford 1995.

\end{document}